\documentclass[12pt]{article}
\usepackage{amsmath,amssymb,amsthm}

\title{The Error in Rayleigh's Approximative Period
}
\author{Mark B. Villarino\\
Depto.\ de Matem\'atica, Universidad de Costa Rica,\\
2060 San Jos\'e, Costa Rica}

\date\

\newtheorem{thm}{Theorem}

\newtheorem{prop}{Proposition}
\newtheorem{ex}{Example}
\numberwithin{equation}{section}

\makeatletter
\def\section{\@startsection{section}{1}{\z@}{-3.5ex plus -1ex minus
			  -.2ex}{2.3ex plus .2ex}{\large\bf}}
\def\subsection{\@startsection{subsection}{2}{\z@}{-3.25ex plus -1ex
			  minus -.2ex}{1.5ex plus .2ex}{\normalsize\bf}}
\renewcommand{\@dotsep}{200} 
\makeatother

\renewcommand{\geq}{\geqslant}  
\renewcommand{\leq}{\leqslant}  

\begin{document}

\maketitle

\begin{abstract}
We  obtain a rigorous \emph{a priori} upper and lower bounds to the exact period of the celebrated \textsc{Rayleigh} stretched-string differential equation.\end{abstract}



\section{Introduction}

In his monumental treatise \emph{The Theory of Sound} (see pages 45-46 \cite{LR}) Lord \textsc{Rayleigh} proposes the following problem:
 \begin{quote}
\emph{A horizontal elastic  wire having length $2L_0$ and spring constant $\sigma$ is stretched to length $2L$.  A body, having a mass $m$ much greater than the mass of the wire, is tied to the center of the wire and is put in motion on a vertical line through the center of the wire.  The only forces on the body are produced by tension in the wire (there is no gravity, no damping).  It is required to study the motion of the body.}
\end{quote}(We have used \textsc{Agnew}'s formulation \cite{ag}, p. 17-18.)  Let $y(t)\equiv y$ be the vertical displacement of the mass $m$ in time $t$, where $y$ is positive, $0$, or negative according as the body is above, or at, or below the wire in equilibrium or neutral position.

By the law of \textsc{Hooke} and the definition of a \emph{spring constant}, $\sigma$, the magnitude of the elastic force, or the \emph{tension}, $T(y)$, exerted by\emph{ each half }of the stretched string is equal to:

\begin{eqnarray*}
T(y)& = &\text{initial tension}+\text{new tension} \\
 &=& \sigma\cdot (L-L_0)+\sigma\cdot( \sqrt{L^2+y^2}-L )\\
 & = &\sigma\cdot( \sqrt{L^2+y^2}-L_0 ) 
\end{eqnarray*} 

\noindent and therefore the \emph{vertical component} of the \emph{total }force, exerted by \emph{both halves}, which is the only component acting to move the mass $m$ is given by

\begin{equation}
\label{ }
\rm{vertical \  component \  of \  T(y)}=-2\sigma\cdot (\sqrt{L^2+y^2}-L_0)\frac{y}{\sqrt{L^2+y^2}},
\end{equation}

\noindent where we have employed our \emph{sign} conventions.

But, by \textsc{Newton}'s second law of motion, that same force, with \emph{opposite} sign, is given by $m\frac{d^2y}{dt^2}$.

Therefore, equating the two forces, and bringing the expression for the elastic force on the same side as the force given by \textsc{Newton}'s second law, we find that the algebraic sum of the forces is zero and dividing by $m$, we conclude that  \emph{the equation of motion} of the mass $m$ is given by:

\begin{equation}
\label{RS}
\boxed{\frac{d^2y}{dt^2}+\Bigl(\frac{2\sigma}{m}\frac{\sqrt{L^2+y^2}-L_0}{\sqrt{L^2+y^2}}\Bigr)y=0}
\end{equation}
This differential equation \eqref{RS} is nonlinear and very complicated, and we have never seen an \emph{exact} treatment of it in the literature.

\textsc{Rayleigh}, himself, simplifies the equation by the following physical reasoning;

\begin{quotation}``The tension of the string in the position of equilibrium depends on the amount of the stretching to which it has been subjected.  In any other position the tension is greater; but we limit ourselves to the case of vibrations so small that the additional stretching is a negligible fraction of the whole.  On this condition, \emph{the tension may be treated as a constant.}"\end{quotation}

(Our italics.)We just saw that the magnitude of the tension, $T(y)$, (or the force) from each half of the string on $m$ is given by
\begin{equation}
\label{T}
T(y)=\sigma(\sqrt{L^2+y^2}-L_0)
\end{equation}provided the stretching is well within elastic limits.  So, \textsc{Rayleigh} \emph{effectively takes }$y=0$ in \eqref{T} so that the equation \eqref{RS} becomes

\begin{equation}
\label{LRS}
\frac{d^2y}{dt^2}+\frac{2\sigma}{m}\frac{L-L_0}{L}y=0,
\end{equation} 

\noindent or, using the notation \eqref{T} with $y=0$, i.e., 

$$T:=T(0)=\sigma (L-L_0)$$ 

\noindent we obtain Rayleigh's approximative differential equation:

\begin{equation}
\label{LRS1}
\boxed{\ddot{y}+\frac{2T}{mL}y=0}
\end{equation}

\noindent where the ``dots" refer to time derivatives.  If we assume that the initial displacement at time $t=0$ is $y_0$, then  the solution to \eqref{LRS1} is 

\begin{equation}
\label{as}
\boxed{y(t)=y_0\cos\left(\sqrt{\frac{2T}{mL}}t\right)}
\end{equation}

\noindent which is \emph{simple harmonic motion} with period $\overline{P}$ given by:
\begin{equation}
\label{aP}
\boxed{\overline{P}=\frac{2\pi}{\sqrt{\frac{2T}{mL}}}}
\end{equation}

\textsc{Rayleigh} does not discuss \emph{how close} his approximate period $\overline{P}$ \eqref{aP} is to the \emph{exact} period $P$.  Indeed, we have been unable to to find any error analysis of \textsc{Rayleigh}'s solution in the literature.  Yet, such a famous and classical differential equation merits an investigation into the accuracy of its approximate solution.

In modern mathematics, any approximation, to be acceptable, requires an upper bound on the error to show how \emph{good} it is, and a lower bound to show how \emph{bad} it is.

Therefore we offer the following theorems to fill this gap.

\begin{thm}
\label{bounds}
Let $P$ be the true period of oscilation of the mass $m$.  Then \textsc{Rayleigh}'s approximate period $\overline{P}$ \textbf{\emph{overestimates}} the true period, $P$.  Indeed, the following inequalities are valid:  If we define the \textbf{relative error} $R$ in the approximation $P\approx\overline{P}$ by
$$
R:=\frac{P-\overline{P}}{P}
$$
then
$$
R<0
$$
and, in fact
\begin{equation}
\label{Bounds}
\boxed{ \sqrt{\frac{\frac{1}{L_0}-\frac{2}{L+\sqrt{L^2+y_0^2}}}{\frac{1}{L_0}-\frac{1}{L}}}<\frac{\overline{P}}{P}<\sqrt{\frac{\frac{1}{L_0}-\frac{1}{\sqrt{L^2+y_0^2}}}{\frac{1}{L_0}-\frac{1}{L}}}}
\end{equation}
\end{thm}

  These inequalities show the surprising result that as long as Hooke's law is applicable, the bound on $R$ only depends on the \emph{ratios} $\frac{L_0}{L}$ and $\frac{y_0}{L}$ and \emph{not} on the mass of the bob nor and Young's modulus of the material.   
  
  We can combine these  radical function bounds into a simple and elegant formula for the (absolute value of) the relative error which displays is true dependence on the various variables in play.
  
  \begin{thm} 
  \label{formula for the relative error}
  \begin{equation}
\label{relative error}
\boxed{|R|=\Theta\cdot\frac{L_0}{L-L_0}\cdot\left(\frac{y_0}{2L}\right)^2}  
  \end{equation}
  where
  \begin{equation}
  \label{bounds on proportion}
  \frac {\frac{1}{2}-\left(\frac{y_0}{2L}\right)^2}{1+\frac{1}{32}\frac{L_0}{L-L_0}}\leq\Theta\leq 1.
 \end{equation}
 So, the (absolute value of) the relative error is \textbf{directly proportional} to the square of the initial fractional displacement $\left(\frac{y_0}{2L}\right)^2$ and \textbf{inversely proportional} to the initial stretch $L-L_0$.  The factor of proportionality $\Theta$ is always greater than $\frac {\frac{1}{2}-\left(\frac{y_0}{2L}\right)^2}{1+\frac{1}{32}\frac{L_0}{L-L_0}}$ but always is less than $1$
.  \end{thm}

Thus, notwithstanding the analytical complexity of the radical function bounds, their behavior is limited by the simple bounds \eqref{bounds on proportion} on $\Theta$.
  
   We believe the theorem and inequalities to be new.

\section {The Formula for the Period}


The differential equation \eqref{RS} does not have the variable $t$ nor the first derivative,  $\dot{y}$ ,appearing explicity.  Therefore, we can transform it into a separable equation for the \emph{velocity}, 

\begin{equation}
\label{v}
v:=\dot{y}
\end{equation}

\noindent since

\begin{equation}
\label{yddot}
\ddot{y}=\frac{dv}{dt}=\frac{dv}{dy}\frac{dy}{dt}=v\frac{dv}{dy}
\end{equation}

Substituting the right-hand side of \eqref{yddot} into \eqref{RS}, transposing the expression for the tension to the right-hand side, separating the variables and forming the indefinite integral of both sides we obtain

\begin{equation}
\label{ }
\frac{v^2}{2}=-\frac{\sigma}{2m}\left(\frac{y^2}{2}-L_0\sqrt{L^2+y^2}\right)+C
\end{equation}

\noindent where $C$ is the constant of integration.  At time $t=0$, the conditions are $y(0)=y_0$ and $v(0)=0$.  Therefore, evaluating the constant $C$, recalling that $v=\frac{dy}{dt}$,  taking the square root of both sides, and using a \emph{negative} sign because the displacement, $y(t)$, is a \emph{decreasing} function during the first quarter oscilation, we obtain:

\begin{equation}
\label{ }
\frac{dy}{dt}=-\sqrt{\frac{2\sigma L_0}{m}(y_0^2-y^2)\left(\frac{1}{L_0}-\frac{2}{\sqrt{L^2+y^2}+\sqrt{L^2+y_0^2}}\right).}
\end{equation}

Multiplying the quarter-period by $4$, we obtain the following formula for the full period, $P$, of oscilation of the mass:

\begin{thm}
The \textbf{true period}, $P$, of oscillation of the mass $m$ is given by:
\begin{equation}
\label{P}
\boxed{P=4\sqrt{\frac{m}{2\sigma L_0}}\int_{0}^{y_0}\frac{1}{\sqrt{y_0^2-y^2}\sqrt{\frac{1}{L_0}-\frac{2}{\sqrt{L^2+y^2}+\sqrt{L^2+y_0^2}}}}~dy.}
\end{equation}
\qed
\end{thm}

This integral is a very complicated elliptic integral and does not permit an elementary exact computation.

Instead, we turn our attention to the upper and lower bounds for the true period, $P$.


\section{Upper and Lower Bounds for the True Period}

First, we find upper and lower bounds for the period $P$:
\begin{prop}
The true period $P$ is subject to the following inequalities:

\begin{equation}
\label{TP}
\boxed{\frac{2\pi\sqrt{\frac{m}{2\sigma L_0}}}{\sqrt{\frac{1}{L_0}-\frac{1}{\sqrt{L^2+y_0^2}}}}\leq P\leq\frac{2\pi\sqrt{\frac{m}{2\sigma L_0}}}{\sqrt{\frac{1}{L_0}-\frac{2}{L+\sqrt{L^2+y_0^2}}}}}
\end{equation}
\end{prop}

\begin{proof}[Proof of the bounds on the true period]
The formula for the true period $P$ is (\eqref{P}):
$$
P=4\sqrt{\frac{m}{2\sigma L_0}}\int_{0}^{y_0}\frac{1}{\sqrt{y_0^2-y^2}\sqrt{\frac{1}{L_0}-\frac{2}{\sqrt{L^2+y^2}+\sqrt{L^2+y_0^2}}}}~dy.
$$
We will find bounds on the second radical in the denominator of the integrand which do not depend on $y$. 
To this end we note that since $0\leq y\leq y_0$ and $0\leq L_0\leq L$, we conclude that

$$
L+\sqrt{L^2+y_0^2}\leq \sqrt{L^2+y^2}+\sqrt{L^2+y_0^2}\leq 2\sqrt{L^2+y_0^2}
$$

\noindent and after some algebraic manipulation we obtain

\begin{equation}
\label{bounds}
\sqrt{\frac{1}{L_0}-\frac{2}{L+\sqrt{L^2+y_0^2}}}\leq\frac{1}{\frac{1}{\sqrt{\frac{1}{L_0}-\frac{2}{\sqrt{L^2+y^2}+\sqrt{L^2+y_0^2}}}}}\leq\sqrt{\frac{1}{L_0}-\frac{1}{\sqrt{L^2+y_0^2}}}
\end{equation}

If we apply this inequality to the denominator in \eqref{P} and use the results
$$
\int_{0}^{y_0}\frac{dy}{\sqrt{y_0^2-y^2}}=\frac{\pi}{2}
$$
we obtain the two bounds in the theorem.
\end{proof}


\section{Upper and lower bounds on the relative error}


Now we use these bounds to obtain the announced bounds on the relative error.
\begin{proof}

Since the relative error
$$
R:=\frac{P-\overline{P}}{P}\equiv 1-\frac{\overline{P}}{P},
$$
it is sufficient to bound the fraction $\frac{\overline{P}}{P}.$

First we reduce the formula for Rayleigh's approximative period:

$$
\overline{P}=\frac{2\pi}{\sqrt{\frac{2T}{mL}}}=2\pi\sqrt{\frac{mL}{2T}}=\frac{2\pi}{\sqrt{\frac{mL}{2\sigma(L-L_0)}}}=2\pi\sqrt{\frac{m}{2\sigma}}\sqrt{\frac{L}{L-L_0}}
$$
and therefore, recalling the exact formula for the period $P$, namely \eqref{P}, we obtain

$$
\frac{\overline{P}}{P}=\frac{2\pi\sqrt{\frac{m}{2\sigma}}\sqrt{\frac{L}{L-L_0}}}{4\sqrt{\frac{m}{2\sigma L_0}}\int_{0}^{y_0}\frac{dy}{\sqrt{y_0^2-y^2}\sqrt{\frac{1}{L_0}-\frac{2}{\sqrt{L^2+y^2}+\sqrt{L^2+y_0^2}}}}.}
$$
which collapses down to  

\begin{equation}
\label{PP1}
\frac{\overline{P}}{P}=\frac{\frac{\pi}{2}\sqrt{\frac{L_0L}{L-L_0}}}{\int_{0}^{y_0}\frac{dy}{\sqrt{y_0^2-y^2}\sqrt{\frac{1}{L_0}-\frac{2}{\sqrt{L^2+y^2}+\sqrt{L^2+y_0^2}}}}.}
\end{equation}

Thus, we must find upper and lower bounds for the denominator in \eqref{PP1}.  But if we combine the bounds in \eqref{TP} with
$$
\frac{LL_0}{L-L_0}=\frac{1}{\frac{1}{L_0}-\frac{1}{L}}
$$
and 
$$
\int_{0}^{y_0}\frac{dy}{\sqrt{y_0^2-y^2}}=\frac{\pi}{2}
$$
and cancel the $\frac{\pi}{2}$ in the numerator and the denominator, we obtain the bounds \eqref{Bounds} in the theorem. 
\end{proof}


\section{Proof of the formula for the relative error.}


First we prove the \emph{upper bound.}

\begin{equation}
\label{upper bound}
|R|\leq\left(\frac{y_0}{2L}\right)^2
\end{equation}

By the right-hand side of \eqref{bounds}
\begin{eqnarray*}
|R| & \leq & \sqrt{\frac{\frac{1}{L_0}-\frac{1}{\sqrt{L^2+y_0^2}}}{\frac{1}{L_0}-\frac{1}{L}}}-1 \\
 & = & \sqrt{\frac{LL_0}{L-L_0}}\cdot\sqrt{\frac{1}{L_0}-\frac{1}{\sqrt{L^2+y_0^2}}}-1\\
 & =&  \sqrt{\frac{LL_0}{L-L_0}}\cdot\sqrt{\frac{1}{L_0}-\frac{1}{L}+\frac{1}{L}-\frac{1}{\sqrt{L^2+y_0^2}}}-1\\
 & = &\sqrt{\frac{LL_0}{L-L_0}}\cdot\sqrt{\frac{1}{L_0}-\frac{1}{L}+\frac{1}{L}-\frac{1}{L}+\frac{y_0^2}{2L^3}-\frac{3y_0^4}{8L^5}+-\cdots}-1\\
 &=&\sqrt{\frac{LL_0}{L-L_0}}\cdot\sqrt{\frac{L-L_0}{L-L_0}+\frac{y_0^2}{2L^3}-\frac{3y_0^4}{8L^5}+-\cdots}-1\\
&=&\sqrt{1+\frac{L_0}{L-L_0}\left(\frac{y_0^2}{2L^2}-\frac{3y_0^4}{8L^4}+-\cdots\right)}-1\\
&<&1+\frac{L_0}{L-L_0}\left(\frac{y_0^2}{4L^2}-\frac{3y_0^4}{16L^4}+-\cdots\right)-1\\
&<&\frac{L_0}{L-L_0}\left(\frac{y_0^2}{4L^2}\right)
\end{eqnarray*}
where we have used the fact that a convergent alternating series whose terms tend monotonically to zero has a sum that is between two consecutive partial sums and that $\sqrt{1+u}<1+\frac{u}{2}.$
This completes the proof of the upper bound.

Now we prove the lower bound.

\begin{equation}
\label{lower bound}
|R|\geq\frac {\frac{1}{2}-\left(\frac{y_0}{2L}\right)^2}{1+\frac{1}{32}\frac{L_0}{L-L_0}}\cdot\frac{L_0}{L-L_0}\left(\frac{y_0}{2L}\right)^2.
\end{equation}

By the right-hand side of \eqref{bounds}
\begin{eqnarray*}
|R| & \geq & \sqrt{\frac{\frac{1}{L_0}-\frac{2}{L+\sqrt{L^2+y_0^2}}}{\frac{1}{L_0}-\frac{1}{L}}}-1 \\
 & = & \sqrt{\frac{LL_0}{L-L_0}}\cdot\sqrt{\frac{1}{L_0}-\frac{2}{L+\sqrt{L^2+y_0^2}}}-1\\
 & =&  \sqrt{\frac{LL_0}{L-L_0}}\cdot\sqrt{\frac{1}{L_0}-\frac{1}{L}+\frac{1}{L}-\frac{2}{L+\sqrt{L^2+y_0^2}}}-1\\
 & = &\sqrt{\frac{LL_0}{L-L_0}}\cdot\sqrt{\frac{1}{L_0}-\frac{1}{L}+\frac{1}{L}-\frac{1}{L}+\frac{y_0^2}{4L^3}-\frac{y_0^4}{8L^5}+-\cdots}-1\\
 &=&\sqrt{\frac{LL_0}{L-L_0}}\cdot\sqrt{\frac{L-L_0}{L-L_0}+\frac{y_0^2}{4L^3}-\frac{y_0^4}{8L^5}+-\cdots}-1\\
&=&\sqrt{1+\frac{L_0}{L-L_0}\left(\frac{y_0^2}{4L^2}-\frac{y_0^4}{8L^4}+-\cdots\right)}-1\\
&>&\sqrt{1+\frac{L_0}{L-L_0}\left(\frac{y_0^2}{4L^2}-\frac{y_0^4}{8L^4}\right)}-1\\
&>&1+\frac{\frac{L_0}{L-L_0}\left(\frac{y_0^2}{4L^2}-\frac{y_0^4}{8L^4}\right)}{2+\frac{\frac{L_0}{L-L_0}\left(\frac{y_0^2}{4L^2}-\frac{y_0^4}{8L^4}\right)}{2}}-1\\
&=&\frac{\frac{1}{2}\frac{L_0}{L-L_0}\left(\frac{y_0^2}{4L^2}-\frac{y_0^4}{8L^4}\right)}{1+\frac{\frac{L_0}{L-L_0}\left(\frac{y_0^2}{4L^2}-\frac{y_0^4}{8L^4}\right)}{4}}\\
&=& \frac{L_0}{L-L_0}\left(\frac{y_0}{2L}\right)^2\frac{\left(\frac{1}{2}-\left(\frac{y_0}{2L}\right)^2\right)}{1+\frac{\frac{L_0}{L-L_0}\left(\frac{y_0^2}{4L^2}-\frac{y_0^4}{8L^4}\right)}{4}}
\end{eqnarray*}

\noindent where we used the inequality $\sqrt{1+u}>1+\frac{u}{2+\frac{u}{2}}.$. Now the function $\frac{y_0^2}{4L^2}-\frac{y_0^4}{8L^4}$ is increasing for $0\leq y_0\leq L$ which means that it takes it maximum value $=\frac{1}{8}$ when $y_0=L.$ Plugging this into the denominator of the last fraction we obtain\begin{equation}
\label{lower bound 1}
|R|>\frac{L_0}{L-L_0}\left(\frac{y_0}{2L}\right)^2\frac{\left(\frac{1}{2}-\left(\frac{y_0}{2L}\right)^2\right)}{1+\frac{1}{32}\frac{L_0}{L-L_0}}\equiv\frac{L_0}{L-\frac{31}{32}L_0}\left(\frac{y_0}{2L}\right)^2\left(\frac{1}{2}-\left(\frac{y_0}{2L}\right)^2\right)
\end{equation}
\noindent The second formula in \eqref{lower bound 1} is an equally elegant form as the original.

This completes the proof of the formula for the relative error.


\section{Examples}


We will finish with some examples.

\begin{ex}
A wire of length $1152$ milimeters is stretched to a length of $1250$ millimeters.  Its spring constant is $1000$  Newtons per meter.  A mass of $0.5$ kilograms is attached to its center.  The mass is displaced vertically $83$ milimeters and released.  It is required to find its \textbf{period of oscillation.}

\emph{The exact formula for the period gives $34.1941$ seconds while Rayleigh's approximative period gives $35.4807$ seconds.  Therefore the absolute value of the relative error is $3.4\%$.  The bounds in \eqref{bounds on proportion}, give
$$
2.5\%\leq 3.4\%\leq 5\%
$$
which verifies the bounds in the theorem in this case.}
\end{ex}

The following example explains an anomaly in the applicability of Rayleigh's approximation.

\begin{ex}
A steel wire one meter long has a tension of $100$ Newtons.  A mass of $0.5$ kilograms is attached to the center and is displaced $2$ centimeters transversely to the wire.  It is required to find the period.  The spring constant of steel is $157,000$ Newtons per meter.
\end{ex}

Rayleigh's approximate period is $0.2221$ seconds while the true period is $0.18308$ seconds.  This gives a relative error of $25\%$!!  This example shows that Rayleigh's approximation can be wildly inaccurate.

Our bounds on the relative error permit us to explain \emph{why} Rayleigh's approximation is so bad in this case.
Namely, for a given displacement, $\frac{y_0}{2L}$,  if \emph{the stretch is tiny}, that is if $L\approx L_0$, the relative error will be very large.  In this particular case the stretch is $\frac{1}{1570}=0.00064$ meters, and thus the lower bound is $15.7\%$ which shows that the true relative error will be even larger.  And, indeed, above we found it to be about $25\%$.

\subsubsection*{Acknowledgment}
 Support from the Vicerrector\'{\i}a de Investigaci\'on of the 
University of Costa Rica is acknowledged.



\begin{thebibliography}{}

\bibitem{ag}
Ralph Palmer Agnew
\textit{Differential Equations}, second edition
McGraw-Hill Inc., New York, 1960.




\bibitem{LR}
John William Strutt
\textit{Theory of Sound, Vol 1},
MacMillan and Co., New York, 1877.





\end{thebibliography}
\end{document}